\documentclass[11pt]{article}
\usepackage[T1]{fontenc}
\usepackage{lmodern,amsmath,amssymb,amsthm,mathtools}
\usepackage[margin=1in]{geometry}
\usepackage{needspace}
\usepackage[hidelinks]{hyperref}
\newtheorem{theorem}{Theorem}[section]
\newtheorem{proposition}[theorem]{Proposition}
\newtheorem{lemma}[theorem]{Lemma}
\newtheorem{corollary}[theorem]{Corollary}
\newcommand{\C}{\mathbb C}
\newcommand{\R}{\mathbb R}
\newcommand{\Ric}{\operatorname{Ric}}
\title{Radial symmetry for the positive K\"ahler--Einstein Dirichlet problem on the ball}
\author{%
Hongyu Chen\thanks{School of Mathematics, Sichuan University,
Chengdu 610064, China.
Email: \texttt{hongyu.chern@gmail.com}.}
\qquad
Li Ma\thanks{Corresponding author.
School of Mathematics and Physics,
University of Science and Technology Beijing,
Beijing 100083, China.
Email: \texttt{nuslma@126.com}.}
}
\date{June 2026}
\begin{document}
\maketitle

\begin{abstract}
We prove radial symmetry for every strictly plurisubharmonic function $u$
on the unit ball in $\C^n$, $n\ge2$, that is smooth up to the boundary,
vanishes there, and has complex Hessian determinant $\mu e^{-u}$ for a
constant $\mu>0$.
The known Fubini--Study family therefore exhausts all such solutions.
For the corresponding complex Monge--Amp\`ere mean field equation, this
establishes the Berman--Berndtsson symmetry conjecture in the smooth
class and gives uniqueness throughout the known existence range.
The classification also yields the sharp Moser--Trudinger inequality
on the ball in Cegrell's finite-energy class.


\end{abstract}

\section{Introduction}
We consider the positive K\"ahler--Einstein equation on the unit ball with
the zero boundary condition. Since the problem is invariant under unitary
transformations, it is natural to ask whether every solution that is smooth
up to the boundary is radial. Actually, Berman and Berndtsson conjectured that every solution of this problem is radial, and
hence coincides with the known radial solution \cite[Conjecture~7.6]{BB2022}. We remark that
this conjecture asks for radial symmetry of arbitrary solutions, rather than
construction or classification within the radial class.

Let $B$ be the unit ball in $\C^n$, where $n\ge2$. In the standard
coordinates $z_j=x_j+iy_j$, write
$\partial_j=\partial_{z_j}=(\partial_{x_j}-i\partial_{y_j})/2$ and
$\partial_{\bar j}=\partial_{\bar z_j}$ for its conjugate.
For derivatives of a smooth function $h$, we write $h_j=\partial_jh$
and $h_{\bar j}=\partial_{\bar j}h$. The complex Hessian of $u$ has
entries $u_{i\bar j}=\partial_i\partial_{\bar j}u$. We consider
\begin{equation}\label{eq:main}
 \det(u_{i\bar j})=\mu e^{-u}\quad\text{in }B,
 \qquad u=0\quad\text{on }\partial B,
\end{equation}
where $\mu>0$ and $u\in C^\infty(\overline B;\R)$ is strictly
plurisubharmonic, meaning that its complex Hessian is positive definite
in $B$. The K\"ahler form $\omega_u=i\partial\bar\partial u$ satisfies
$\Ric(\omega_u)=\omega_u$, since its Ricci form is
$-i\partial\bar\partial\log\det(u_{i\bar j})$. We use $dd^c=(i/2\pi)\partial\bar\partial$ and the probability measure
$dV_0=(n!/\pi^n)dV$ on $B$, where $dV$ is Euclidean volume. For
smooth $u$, $(dd^cu)^n=\det(u_{i\bar j})dV_0$ and
$\omega_u=2\pi dd^cu$.
Our aim here is to show the radial symmetry of any solution to the problem
\eqref{eq:main}.

In complex dimension one, radial symmetry follows from moving
planes \cite{GNN1979}; see \cite[Section~7.3]{BB2022}.
For $n\ge2$, general real hyperplane reflections do not preserve
the complex Monge--Amp\`ere operator. The smooth Dirichlet theory
for prescribed positive densities \cite{CKNS1985} does not give
uniqueness for the decreasing density in \eqref{eq:main}; indeed,
two radial solutions can occur for the same $\mu$.

The explicit radial solutions of \eqref{eq:main} are the Fubini--Study
potentials displayed in \eqref{eq:radial} below and recorded by Li
\cite{Li2016}.
 For bounded plurisubharmonic potentials we use the
Bedford--Taylor Monge--Amp\`ere measure. The mean field equation with
coupling parameter $\gamma>0$ is
\begin{equation}\label{eq:mean-field}
 (dd^c\varphi)^n
 =\frac{e^{-\gamma\varphi}}{\int_Be^{-\gamma\varphi}\,dV_0}\,dV_0,
 \qquad \varphi|_{\partial B}=0,
\end{equation}
Its total Monge--Amp\`ere mass is one.
Rescaling Euclidean volume preserves the normalized measure in
\eqref{eq:mean-field}, which is equation~(7.8) of \cite{BB2022} in the
same $dd^c$ convention.

The variational formulation of \eqref{eq:mean-field} incorporates
zero boundary data by approximation. Let $\mathcal E_0(B)$
consist of the bounded plurisubharmonic functions
$v$ on $B$ with $\lim_{z\to\zeta}v(z)=0$ at every $\zeta\in\partial B$
and $\int_B(dd^cv)^n<\infty$. Cegrell's finite-energy class
$\mathcal E^1(B)$ consists of the plurisubharmonic functions $v$ on $B$
that admit a sequence $v_j\in\mathcal E_0(B)$ decreasing pointwise to
$v$ with
$$
 \sup_{j\ge1}\int_B(-v_j)(dd^cv_j)^n<\infty.
$$
For $v\in\mathcal E^1(B)$, the measure $(dd^cv)^n$ is the weak limit
on $B$ of $(dd^cv_j)^n$, tested against compactly supported continuous
functions. The limit is independent of the defining sequence, and
$\int_B(-v)(dd^cv)^n$ is finite
\cite[Section~2 and Theorems~3.5 and~3.8]{Cegrell1998}.

For $\gamma>0$ and $\varphi\in\mathcal E^1(B)$, the functional of
\cite[Section~4.2]{BB2022} in our volume normalization is
$$
 G_\gamma(\varphi)
 =\frac1\gamma\log\int_B e^{-\gamma\varphi}\,dV_0
   -\frac1{n+1}\int_B(-\varphi)(dd^c\varphi)^n.
$$
Replacing $dV$ by $dV_0$ adds a constant independent of $\varphi$.
A bounded domain is hyperconvex if it admits a negative continuous
plurisubharmonic exhaustion. In the subcritical range $0<\gamma<n+1$,
Berman and Berndtsson showed that the corresponding functional on
any bounded hyperconvex domain attains its maximum at a continuous
solution of the mean field equation \cite[Theorem~7.2]{BB2022}.
They record smooth attainment for domains with smooth strictly
pseudoconvex boundary \cite[Remark~7.3]{BB2022}, using the iteration
of Guedj, Kolev and Yeganefar \cite{GKY2013}. On $B$, the function
$|z|^2-1$ gives both an exhaustion and a smooth strictly
plurisubharmonic defining function.
On the ball, the explicit radial family
already gives smooth existence throughout this interval. Li proved
that there is no smooth solution on the ball when $\gamma\ge n+1$,
without assuming symmetry, using a complex Pohozaev identity
\cite{Li2016}. Thus the sharp existence interval
$0<\gamma<n+1$ was known before the present work. For each
parameter in this interval, the radial solution was also known to
be unique. See \cite[Section~7.3]{BB2022} and \cite{Li2016}.

Two earlier uniqueness results are particularly relevant.
Guedj, Kolev and Yeganefar considered smooth bounded strictly
pseudoconvex domains containing the origin and invariant under the
diagonal circle action $z\mapsto e^{i\theta}z$, $\theta\in\R$. They proved uniqueness
among $S^1$-invariant solutions when the domain is strictly convex for
the metric of one such solution \cite[Theorem~25]{GKY2013}. The lower
bound on the first Dirichlet Laplace eigenvalue used in their
uniqueness argument fails on sufficiently large Fubini--Study balls
\cite{GKYspectral2013}. Lu and Phung proved uniqueness for
sufficiently small positive coupling parameters on bounded
hyperconvex domains, for probability densities in $L^p$ with
$p>1$ \cite[Theorem~1.2]{LP2025}. Their result applies to a
broader class of domains, densities and bounded solutions, but
does not give uniqueness on the whole subcritical interval on
the ball.

The following theorem shows that the known radial family exhausts
all solutions of \eqref{eq:main} that are smooth up to the boundary.

\begin{theorem}\label{thm:main}
Let $n\ge2$, let $B\subset\C^n$ be the unit ball, and let $\mu>0$.
If $u\in C^\infty(\overline B;\R)$ has positive definite complex Hessian
in $B$ and satisfies
$$
 \det(u_{i\bar j})=\mu e^{-u}\quad\text{in }B,
 \qquad u=0\quad\text{on }\partial B,
$$
then $u$ is radial and, for some $b>0$,
\begin{equation}\label{eq:radial}
 u(z)=(n+1)\log\frac{1+b|z|^2}{1+b},
 \qquad
 \mu=\frac{(n+1)^n b^n}{(1+b)^{n+1}}.
\end{equation}
Thus there are exactly two such solutions when
$0<\mu<n^n/(n+1)$, exactly one when $\mu=n^n/(n+1)$, and none when
$\mu>n^n/(n+1)$.
\end{theorem}

We remark that for $0<\mu<n^n/(n+1)$, the relation
$\mu=\frac{(n+1)^n b^n}{(1+b)^{n+1}}$ holds for exactly two positive
values of $b$, as shown in Lemma~\ref{lem:positive-roots}. The relation with  \cite[Conjecture~7.6]{BB2022} is exact within the class of
solutions smooth up to the boundary. If $\varphi$ solves
\eqref{eq:mean-field}, then $u=\gamma\varphi$ solves
\eqref{eq:main} with
$\mu=\gamma^n/\int_Be^{-u}\,dV_0$. Conversely, for any solution
$u$ of \eqref{eq:main}, take
$\gamma=(\int_B(dd^cu)^n)^{1/n}$ and $\varphi=u/\gamma$;
then \eqref{eq:mean-field} holds. These transformations preserve
radial symmetry. Thus the symmetry assertion of
Theorem~\ref{thm:main} is equivalent to the smooth case of the
conjecture, and Corollary~\ref{cor:mean-field} states its
consequence in the original mean field formulation.

In the variables of \eqref{eq:mean-field}, the classification gives
uniqueness throughout the subcritical range.

\begin{corollary}\label{cor:mean-field}
Let $B\subset\C^n$ be the unit ball, with $n\ge2$, and let
$\gamma>0$. Consider strictly
plurisubharmonic functions
$\varphi\in C^\infty(\overline B;\R)$ satisfying
$$
 (dd^c\varphi)^n
 =\frac{e^{-\gamma\varphi}}{\int_Be^{-\gamma\varphi}\,dV_0}\,dV_0,
 \qquad \varphi|_{\partial B}=0.
$$
Such a function exists exactly when $0<\gamma<n+1$, and is then unique:
$$
 \varphi(z)=\frac{n+1}{\gamma}
   \log\left(1-\frac\gamma{n+1}
                  +\frac\gamma{n+1}|z|^2\right).
$$
\end{corollary}

\Needspace{12\baselineskip}

For smooth $S^1$-invariant plurisubharmonic potentials on the ball
with zero boundary values, Berman and Berndtsson proved the
Moser--Trudinger inequality with the sharp energy coefficient
$(n+1)^{-(n+1)}$. For all smooth plurisubharmonic potentials with
the same boundary data, their estimate holds with any strictly
larger coefficient \cite[Theorems~1.4 and~1.5]{BB2022}. They
showed how radial symmetry of the subcritical maximizers of
$G_\gamma$ yields the inequality at the critical parameter
$\gamma=n+1$ \cite[Section~7.3]{BB2022}. Applying
Corollary~\ref{cor:mean-field} to these maximizers gives the
following result.

\begin{corollary}\label{cor:moser-trudinger}
Let $B$ be the unit ball in $\C^n$, $n\ge2$. With the above
normalizations of $dd^c$ and $dV_0$, every $u\in\mathcal E^1(B)$ satisfies
$$
 \log\int_B e^{-u}\,dV_0
 \le \frac{1}{(n+1)^{n+1}}\int_B(-u)(dd^cu)^n
     +\sum_{j=1}^n\frac1j.
$$
The energy coefficient is optimal among inequalities with a finite
additive constant, and $\sum_{j=1}^n 1/j$ is the best additive
constant for that coefficient.
\end{corollary}

Delano\"e studied radial solutions to real and complex Monge--Amp\`ere
boundary-value problems in balls \cite{Delanoe1985}. For planar mean
field equations on annuli, Ma obtained comparison under Kelvin
inversion and a radial-derivative dichotomy under sign and
inversion monotonicity assumptions on the coefficient
\cite[Theorem~3]{Ma2009}.

Schwarz symmetrization preserves plurisubharmonicity for $S^1$-invariant
potentials on balanced domains \cite[Theorem~2.3]{BB2014}.
On $B$, it also does not increase Monge--Amp\`ere energy for
such potentials in $C(\overline B)$ vanishing on $\partial B$
\cite[Theorem~2.4]{BB2014}. Extensions concern fiberwise circle-invariant
potentials on disk bundles over Fano manifolds under curvature
hypotheses \cite{Wu2020}, and quasi-circular invariant
potentials \cite{DHZ2026}.

Exponential-integrability estimates in finite-energy classes
and on pseudoconvex domains were obtained in
\cite{AC2019,WWZ2020,DGL2021}; the local estimate in
\cite{DGL2021} uses a strictly subcritical exponential
coefficient. For more general decreasing densities, flow
methods give variational existence under sublinear or
superlinear growth assumptions \cite{CLMZ2026}.
Badiane and Zeriahi proved uniqueness under a uniform bound
on the negative slope of the $n$-th root of the density in
terms of the first Dirichlet eigenvalue of the complex Monge--Amp\`ere operator
\cite[Proposition~4.2]{BZ2023}; the corrected gradient estimate
appears in \cite{BLZ2026}. The exponential density in
\eqref{eq:main} does not satisfy that slope bound on the whole
negative half-line. These results do not supply the
classification asserted in Theorem~\ref{thm:main}.

The proof of Theorem~\ref{thm:main} proceeds as follows.
In Section~2, we derive the mass and Pohozaev identities and
the Green formula in Lemma~\ref{lem:boundary}. Comparing with
the radial solution of smaller boundary slope gives a lower
bound for the outward normal derivative in
Lemma~\ref{lem:lower-slope}. In Section~3, we study a logarithmic auxiliary function defined
using the gradient norm in the solution metric (see
\cite{Enache2014} for auxiliary functions using the Euclidean
gradient in real Monge--Amp\`ere equations).
Proposition~\ref{prop:P} establishes its subharmonicity with
respect to the solution metric, and
Proposition~\ref{prop:log-flux} computes its boundary flux.
In Section~\ref{sec:normalized}, we use the
comparison bound and the Pohozaev identity to show that the
flux vanishes. The resulting holomorphic gradient yields the
radial classification in Theorem~\ref{thm:main}. We then prove
the mean field classification and the sharp Moser--Trudinger
inequality in Corollaries~\ref{cor:mean-field}
and~\ref{cor:moser-trudinger}.

\section{Boundary identities and comparison}

Fix a solution of \eqref{eq:main}. Write $g=(g_{i\bar j})$ for its
K\"ahler metric, with $g_{i\bar j}=u_{i\bar j}$, and let
$g^{i\bar j}$ be the $(j,i)$ entry of the inverse matrix.
We use the complex Laplacian
$\Delta_gh=\sum_{i,j}g^{i\bar j}h_{i\bar j}$ and, for real $h$,
$|\partial h|_g^2=\sum_{i,j}g^{i\bar j}h_i h_{\bar j}$.
We write $\langle\cdot,\cdot\rangle_g$ and $|\cdot|_g$ for the
Hermitian tensor pairings and norms induced by $g$.
Covariant derivatives use the Levi--Civita connection of $g$, and
$\nabla^{2,0}h$ is the $(2,0)$ part of the covariant Hessian of $h$.
Thus $\Delta_gu=n$ and $\Ric(\omega_u)=\omega_u$.
Let $\nu(z)=z$ be the Euclidean outward unit normal on $\partial B$
and write $h_\nu$ for the derivative of $h$ in that direction.
The subharmonic maximum principle and Hopf lemma give $u<0$ in $B$
and $u_\nu>0$ on $\partial B$. Continuity and $\det g=\mu$ on the boundary imply
that $g$ is positive definite on $\overline B$, even if positivity
was initially assumed only in $B$.

Let $J$ be the standard complex structure on $\C^n$ and let
$H\subset T\partial B$ be the Euclidean orthogonal complement of
$J\nu$ in $T\partial B$. We use the Euclidean inner product
$\langle\cdot,\cdot\rangle$ on the boundary. Write $\nabla_H$ for the projection of the
spherical gradient onto $H$ and $\Delta_H$ for the trace of the
spherical Hessian over $H$. Both operators act on boundary traces
when applied to functions defined on $\overline B$. Let $d\sigma_0$
be the probability surface measure on $\partial B$. The field $J\nu$
is divergence-free and its integral curves are unit-speed great
circles. Writing $\Delta_{\partial B}$ for the spherical Laplacian,
we have $\Delta_H=\Delta_{\partial B}-(J\nu)^2$ and, for real
$h,k\in C^\infty(\partial B)$,
\begin{equation}\label{eq:horizontal-ibp}
 \int_{\partial B}h\Delta_Hk\,d\sigma_0
 =-\int_{\partial B}\langle\nabla_Hh,\nabla_Hk\rangle\,d\sigma_0.
\end{equation}

The mass identity follows from the cofactor divergence formula.
The Pohozaev identity below is the specialization of Li's complex
Pohozaev identity to \eqref{eq:main} on the unit ball \cite{Li2016}.
We include its derivation in our volume conventions, together with
the Green formula required in Proposition~\ref{prop:log-flux}.

\begin{lemma}\label{lem:boundary}
Let $u\in C^\infty(\overline B;\R)$ be strictly plurisubharmonic
and solve \eqref{eq:main}. Then
\begin{align}
 &\int_B\det(u_{i\bar j})\,dV_0
 =\int_{\partial B}\left(\frac{u_\nu}{2}\right)^n d\sigma_0,\notag\\
 &\int_{\partial B}\left(\frac{u_\nu}{2}\right)^n
       \left(1-\frac{u_\nu}{2(n+1)}\right)d\sigma_0
 =\mu.\label{eq:mass-pohozaev}
\end{align}
For every real $h\in C^\infty(\overline B)$,
\begin{equation}\label{eq:flux}
 \int_B\Delta_gh\,\det(u_{i\bar j})\,dV_0
 =\frac n2\int_{\partial B}\left(\frac{u_\nu}{2}\right)^{n-1}
 \left(h_\nu-
 \frac{\langle\nabla_Hu_\nu,\nabla_Hh\rangle}{2u_\nu}\right)d\sigma_0.
\end{equation}
\end{lemma}

\begin{proof}
At a boundary point, choose unitary coordinates in which the point
is $e_n=(0,\ldots,0,1)$. For tangent vectors $X,Y$ to the sphere, the boundary
condition gives
$$
 D^2u(X,Y)=u_\nu\langle X,Y\rangle,\qquad
 D^2u(X,\nu)=Xu_\nu,
$$
where $D^2u$ is the Euclidean Hessian. Writing
$X_\alpha=\partial_{x_\alpha}$ and
$Y_\alpha=\partial_{y_\alpha}$, we have, for $\alpha,\beta<n$,
$$
 g_{\alpha\bar\beta}=\frac{u_\nu}{2}\delta_{\alpha\beta},\qquad
 g_{\alpha\bar n}=
 \frac{X_\alpha u_\nu-iY_\alpha u_\nu}{4}.
$$
The block inverse formula and $\det g=\mu$ give
$$
 (\det g)g^{n\bar n}=\left(\frac{u_\nu}{2}\right)^{n-1},\qquad
 (\det g)g^{n\bar\alpha}=
 -\frac14\left(\frac{u_\nu}{2}\right)^{n-2}
 (X_\alpha u_\nu-iY_\alpha u_\nu).
$$
The cofactor identities
$\sum_i\partial_i((\det g)g^{i\bar j})=0$ follow from the symmetry
of the third derivatives of $u$. In the stated probability
normalizations, the divergence theorem therefore gives
$$
 \int_B\Delta_gh\,\det g\,dV_0
 =n\int_{\partial B}\operatorname{Re}
 \left(\sum_{i,j}\bar z_i(\det g)g^{i\bar j}h_{\bar j}\right)d\sigma_0.
$$
At $e_n$, the boundary integrand satisfies
$$
 \operatorname{Re}\sum_j(\det g)g^{n\bar j}h_{\bar j}
 =\frac12\left(\frac{u_\nu}{2}\right)^{n-1}
 \left(h_\nu-
 \frac{\langle\nabla_Hu_\nu,\nabla_Hh\rangle}{2u_\nu}\right),
$$
which proves \eqref{eq:flux}. For real
$v,h\in C^\infty(\overline B)$, integrate $v\Delta_gh$ and
$h\Delta_gv$ against $(\det g)dV_0$ by parts and subtract.
The real parts of the interior gradient terms cancel, giving
\begin{align*}
 &\int_B(v\Delta_gh-h\Delta_gv)\det g\,dV_0=\frac n2\int_{\partial B}\left(\frac{u_\nu}{2}\right)^{n-1}
 \left(vh_\nu-hv_\nu-
 \frac{\langle\nabla_Hu_\nu,v\nabla_Hh-h\nabla_Hv\rangle}
 {2u_\nu}\right)d\sigma_0.
\end{align*}
Taking $h=u$ in \eqref{eq:flux} and using $\nabla_Hu=0$ and
$\Delta_gu=n$ proves the mass identity in
\eqref{eq:mass-pohozaev}.

For the Pohozaev identity, the Green identity with $v=u$ yields,
for every real smooth variation $h$, with no restriction on its
boundary values,
\begin{align*}
 \left.\frac d{d\varepsilon}\right|_{\varepsilon=0}
 \frac1{n+1}\int_B(u+\varepsilon h)
       \det((u+\varepsilon h)_{i\bar j})\,dV_0&=\frac1{n+1}\int_B(h+u\Delta_gh)\det g\,dV_0\\
 &\quad=\int_Bh\det g\,dV_0-
 \frac n{n+1}\int_{\partial B}h\left(\frac{u_\nu}{2}\right)^n d\sigma_0.
\end{align*}
Set $u_s(z)=u(sz)$ for $0<s\le1$. A change of variables gives
$$
 \frac1{n+1}\int_Bu_s\det((u_s)_{i\bar j})\,dV_0
 =\frac1{n+1}\int_{sB}u\det g\,dV_0.
$$
The derivative at $s=1$ from below is zero because $u=0$ on
$\partial B$. Since $\left.\partial_su_s\right|_{s=1}
=\sum_i(z_iu_i+\bar z_iu_{\bar i})$ and this equals $u_\nu$ on
the sphere, the variation formula gives
$$
 \int_B\sum_i(z_iu_i+\bar z_iu_{\bar i})\det g\,dV_0
 =\frac{2n}{n+1}\int_{\partial B}
       \left(\frac{u_\nu}{2}\right)^{n+1}d\sigma_0.
$$
The real Euler vector field
$\sum_i(z_i\partial_{z_i}+\bar z_i\partial_{\bar z_i})$ has Euclidean
divergence $2n$. Applying this field to $\det g=\mu e^{-u}$ and using
the divergence theorem gives
$$
 \int_B\sum_i(z_iu_i+\bar z_iu_{\bar i})\det g\,dV_0
 =2n\left(\int_B\det g\,dV_0-\mu\right).
$$
Combining these two identities with the mass formula proves
the second identity in \eqref{eq:mass-pohozaev}.
\end{proof}

\begin{lemma}\label{lem:positive-roots}
Let $n\ge2$. The function $s\mapsto s^n(1-s/(n+1))$ on $(0,\infty)$
attains its maximum $n^n/(n+1)$ only at $s=n$. For
$0<\mu<n^n/(n+1)$, the equation $s^n(1-s/(n+1))=\mu$ has exactly two
positive roots, one in $(0,n)$ and one in $(n,n+1)$. In the same range
of $\mu$, the parameter relation in \eqref{eq:radial} has exactly two
positive roots in $b$.
\end{lemma}

\begin{proof}
The derivative of $s^n(1-s/(n+1))$
is $s^{n-1}(n-s)$. The function is therefore strictly increasing on
$(0,n)$ and strictly decreasing on $(n,\infty)$, with maximum
$n^n/(n+1)$ at $s=n$. Since the function vanishes at $s=0$ and $s=n+1$,
the intermediate value theorem and strict monotonicity give the claimed
root count in $s$.

The substitution $s=(n+1)b/(1+b)$ is a strictly increasing bijection
from $(0,\infty)$ onto $(0,n+1)$ and gives
$$
 s^n\left(1-\frac{s}{n+1}\right)
 =\frac{(n+1)^n b^n}{(1+b)^{n+1}}.
$$
Both roots in $s$ lie in $(0,n+1)$, so the parameter relation in
\eqref{eq:radial} has exactly two positive roots in $b$.
\end{proof}

Applying the maximum bound in Lemma~\ref{lem:positive-roots} to
\eqref{eq:mass-pohozaev} gives $\mu\le n^n/(n+1)$. If equality holds,
the continuous integrand equals its maximum everywhere on $\partial B$
because $d\sigma_0$ is a probability measure with full support. Thus
$u_\nu/2\equiv n$ on $\partial B$.

\begin{lemma}\label{lem:lower-slope}
Let $0<\mu<n^n/(n+1)$, and let $a\in(0,n)$ be the smaller
positive root of $s^n(1-s/(n+1))=\mu$.
If $u\in C^\infty(\overline B;\R)$ is strictly plurisubharmonic
and solves \eqref{eq:main}, then $u_\nu/2\ge a$ on $\partial B$.
\end{lemma}

\begin{proof}
Use the known Fubini--Study solution \cite{Li2016}
$$
 v(z)=(n+1)\log\left(1-\frac a{n+1}+\frac a{n+1}|z|^2\right).
$$
Its tangential and radial complex Hessian eigenvalues are,
respectively,
$$
 \frac{a}{1-a/(n+1)+a|z|^2/(n+1)},\qquad
 \frac{a(1-a/(n+1))}{(1-a/(n+1)+a|z|^2/(n+1))^2}.
$$
Consequently, $v$ solves \eqref{eq:main} with coefficient
$a^n(1-a/(n+1))=\mu$ and has $v_\nu/2=a$ on the sphere.

Write $(g_v)_{i\bar j}=v_{i\bar j}$ and use
$\Delta_{g_v}h=\sum_{i,j}g_v^{i\bar j}h_{i\bar j}$, with the same
inverse-matrix convention as for $g$. Differentiating the equation
satisfied by $v(e^sz)-2ns$ at $s=0$ gives
$$
 (\Delta_{g_v}+1)\left(\sum_i(z_iv_i+\bar z_iv_{\bar i})-2n\right)=0.
$$
The explicit formula for $v$ shows that
$$
 2n-\sum_i(z_iv_i+\bar z_iv_{\bar i})
 =\frac{2\{n(n+1-a)-a|z|^2\}}{n+1-a+a|z|^2}>0
 \quad\text{on }\overline B,
$$
since $a<n$. Concavity of $\log\det$ gives
$-(u-v)\le\Delta_{g_v}(u-v)$. Dividing $u-v$ by the positive
function in the last display and using its equation yields
\begin{align*}
 &\Delta_{g_v}\left[
 \frac{(n+1-a+a|z|^2)(u-v)}{2\{n(n+1-a)-a|z|^2\}}\right]\\
 &\quad+2\operatorname{Re}\sum_{i,j}g_v^{i\bar j}
 \left(\partial_i\log\frac{2\{n(n+1-a)-a|z|^2\}}
 {n+1-a+a|z|^2}\right)
 \partial_{\bar j}\left[
 \frac{(n+1-a+a|z|^2)(u-v)}{2\{n(n+1-a)-a|z|^2\}}\right]\ge0.
\end{align*}
The quotient is zero on $\partial B$. The coefficients are smooth
on $\overline B$, and the operator is uniformly elliptic with no
zeroth-order term. The maximum principle makes the quotient
nonpositive, hence $u\le v$. Since $v-u$ is nonnegative in $B$
and vanishes on $\partial B$, its outward normal derivative is
nonpositive. Therefore $u_\nu/2\ge a$.
\end{proof}

\section{The logarithmic auxiliary function}
  We now use K\"ahler-Bochner formula to derive the following result, which play an important role in the proof of Theorem \ref{thm:main}.
\begin{proposition}\label{prop:P}
Let $u\in C^\infty(\overline B;\R)$ be strictly plurisubharmonic
and solve \eqref{eq:main}. Then
\begin{enumerate}
\item[(1)] $\displaystyle
 e^{-u/(n+1)}\left(1+\frac{|\partial u|_g^2}{n+1}\right)>1$
on $\overline B$.

\item[(2)] $\displaystyle
 \Delta_g\log\left[e^{-u/(n+1)}
 \left(1+\frac{|\partial u|_g^2}{n+1}\right)-1\right]\ge0$
in $B$.

\item[(3)] If
$
 \int_B\Delta_g\log\left[
 e^{-u/(n+1)}\left(1+\frac{|\partial u|_g^2}{n+1}\right)-1
 \right]\det g\,dV_0=0,
$
then
$ \nabla^{2,0}u=\frac1{n+1}\partial u\otimes\partial u$ in $B$.
\end{enumerate}
\end{proposition}

\begin{proof}
Since $u<0$ in $B$, we have $e^{-u/(n+1)}>1$ there. On
$\partial B$, $u=0$ and $u_\nu>0$, so $|\partial u|_g^2>0$.
Consequently, the function
$$
 e^{-u/(n+1)}\left(1+\frac{|\partial u|_g^2}{n+1}\right)-1
$$
is positive on $\overline B$. By continuity and compactness, it
has a positive minimum there, so its logarithm is smooth on
$\overline B$.

The K\"ahler Bochner identity, together with $\Delta_gu=n$ and
$\Ric(\omega_u)=\omega_u$, gives
$$
 \Delta_g|\partial u|_g^2=n+|\partial u|_g^2+|\nabla^{2,0}u|_g^2.
$$
Here the mixed Hessian and Ricci terms are $n$ and
$|\partial u|_g^2$, respectively, while the derivatives of
$\Delta_gu$ vanish. In holomorphic coordinates normal for $g$
at a fixed point,
\begin{align*}
 \partial_j|\partial u|_g^2&=u_j+\sum_i(\nabla^{2,0}u)_{ij}u_{\bar i},\qquad
 \partial_je^{-u/(n+1)}=-\frac{e^{-u/(n+1)}}{n+1}u_j,\\
 \Delta_ge^{-u/(n+1)}&=e^{-u/(n+1)}
 \left(-\frac n{n+1}+\frac{|\partial u|_g^2}{(n+1)^2}\right).
\end{align*}
The product rule gives
{\small
\begin{align}
 \Delta_g\left[e^{-u/(n+1)}\left(1+\frac{|\partial u|_g^2}{n+1}\right)\right]
 &=e^{-u/(n+1)}\left\{\frac{2|\partial u|_g^2}{(n+1)^2}
 +\frac{|\nabla^{2,0}u|_g^2}{n+1}
 +\frac{|\partial u|_g^4}{(n+1)^3}-\frac2{(n+1)^2}\operatorname{Re}
 \langle\partial u,\partial|\partial u|_g^2\rangle_g\right\}\notag\\
 &=\frac{e^{-u/(n+1)}}{n+1}\left|\nabla^{2,0}u-\frac1{n+1}\partial u\otimes\partial u\right|_g^2.\label{eq:P}
\end{align}}
The completed square uses
$ \operatorname{Re}\langle\partial u,\partial|\partial u|_g^2\rangle_g
 =|\partial u|_g^2+\operatorname{Re}
 \langle\nabla^{2,0}u,\partial u\otimes\partial u\rangle_g.$
In the same normal coordinates,
$$
 \partial_j\left[e^{-u/(n+1)}\left(1+\frac{|\partial u|_g^2}{n+1}\right)\right]
 =\frac{e^{-u/(n+1)}}{n+1}\sum_i
 \left((\nabla^{2,0}u)_{ij}-\frac{u_i u_j}{n+1}\right)u_{\bar i}.
$$
Cauchy--Schwarz and \eqref{eq:P} therefore yield
\begin{equation}\label{eq:kato}
 \left|\partial\left[e^{-u/(n+1)}\left(1+\frac{|\partial u|_g^2}{n+1}\right)\right]\right|_g^2
 \le\frac{e^{-u/(n+1)}|\partial u|_g^2}{n+1}
       \Delta_g\left[e^{-u/(n+1)}\left(1+\frac{|\partial u|_g^2}{n+1}\right)\right].
\end{equation}
The chain rule, \eqref{eq:P} and \eqref{eq:kato} give
{\small
\begin{align}
 &\Delta_g\log\left[e^{-u/(n+1)}\left(1+\frac{|\partial u|_g^2}{n+1}\right)-1\right]=\frac{\displaystyle\Delta_g\left[e^{-u/(n+1)}\left(1+\frac{|\partial u|_g^2}{n+1}\right)\right]}
 {\displaystyle e^{-u/(n+1)}\left(1+\frac{|\partial u|_g^2}{n+1}\right)-1}-\frac{\displaystyle
       \left|\partial\left[e^{-u/(n+1)}\left(1+\frac{|\partial u|_g^2}{n+1}\right)\right]\right|_g^2}
 {\displaystyle\left[e^{-u/(n+1)}\left(1+\frac{|\partial u|_g^2}{n+1}\right)-1\right]^2}\notag\\
 &\quad \quad\ge\frac{e^{-u/(n+1)}-1}{\displaystyle\left[e^{-u/(n+1)}\left(1+\frac{|\partial u|_g^2}{n+1}\right)-1\right]^2}
       \Delta_g\left[e^{-u/(n+1)}\left(1+\frac{|\partial u|_g^2}{n+1}\right)\right]\ge0.\label{eq:log-positive}
\end{align}}

For assertion~(3), the integrand is continuous and nonnegative.
Since $\det g>0$, the vanishing integral implies that the Laplacian
of the logarithm vanishes in $B$. As $e^{-u/(n+1)}>1$ there,
\eqref{eq:log-positive} forces the Laplacian in \eqref{eq:P}
to vanish. Then we conclude that
$\nabla^{2,0}u=\frac1{n+1}\partial u\otimes\partial u$ in $B$.
\end{proof}

\begin{proposition}\label{prop:log-flux}
Let $u\in C^\infty(\overline B;\R)$ be strictly plurisubharmonic
and solve \eqref{eq:main}. Then
\begin{align}
 &\int_B\Delta_g\log\left[e^{-u/(n+1)}\left(1+\frac{|\partial u|_g^2}{n+1}\right)-1\right]\det g\,dV_0\notag\\
 &\quad=-n^2\left\{
 \int_{\partial B}\frac{
 \left(\frac{u_\nu}{2}\right)^n\left(1-\frac{u_\nu}{2(n+1)}\right)-\mu}
 {u_\nu/2}\,d\sigma_0
 +\frac1{16}\int_{\partial B}
 \left(\frac{u_\nu}{2}\right)^{n-3}|\nabla_Hu_\nu|^2\,d\sigma_0
 \right\}.\label{eq:log-flux}
\end{align}
\end{proposition}

\begin{proof}
By Proposition~\ref{prop:P}, the logarithm in \eqref{eq:log-flux}
is smooth on $\overline B$. We compute its boundary flux using
\eqref{eq:flux}, in the unitary coordinates of
Lemma~\ref{lem:boundary}. Write $u_{\nu\nu}=D^2u(\nu,\nu)$.
The block inverse formula and $D^2u(J\nu,J\nu)=u_\nu$ give
\begin{align}
 u_{n\bar n}
 =\frac{\mu}{(u_\nu/2)^{n-1}}
   +&\frac{|\nabla_Hu_\nu|^2}{8u_\nu},\quad
 u_{\nu\nu}
 =\frac{4\mu}{(u_\nu/2)^{n-1}}
   +\frac{|\nabla_Hu_\nu|^2}{2u_\nu}-u_\nu,\notag\\
 |\partial u|_g^2
 &=\left(\frac{u_\nu}{2}\right)^2g^{n\bar n}
   =\frac{(u_\nu/2)^{n+1}}{\mu}.\label{eq:schur}
\end{align}
Keep the coordinate frame fixed along the ray $re_n$.
At $e_n$, $u_\alpha=0$ for $\alpha<n$ and $u_n=u_\nu/2$.
The differentiated gradient factors satisfy
$$
 \partial_\nu u_\alpha=2g_{\alpha\bar n}\quad(\alpha<n),\qquad
 \partial_\nu u_n=\frac12u_{\nu\nu}-\frac i2(J\nu)u_\nu,
$$
and the inverse matrix gives
$$
 \sum_jg^{\alpha\bar j}u_{\bar j}
 =-g_{n\bar\alpha}g^{n\bar n}\quad(\alpha<n),\qquad
 \sum_jg^{n\bar j}u_{\bar j}=
 \frac{u_\nu}{2}g^{n\bar n}.
$$
The term containing $(J\nu)u_\nu$ has zero real part, and
$\sum_{\alpha<n}|g_{\alpha\bar n}|^2
=|\nabla_Hu_\nu|^2/16$. Hence the two differentiated gradient
factors contribute
$$
 2\operatorname{Re}\sum_{i,j}
 g^{i\bar j}(\partial_\nu u_i)u_{\bar j}
 =\left(\frac{u_\nu u_{\nu\nu}}2
       -\frac{|\nabla_Hu_\nu|^2}{4}\right)g^{n\bar n}.
$$
The remaining term differentiates the inverse metric. Since the
gradient has only a normal component at $e_n$, we obtain
\begin{equation}\label{eq:qnu-first}
 \partial_\nu|\partial u|_g^2
 =\left(\frac{u_\nu u_{\nu\nu}}2
       -\frac{|\nabla_Hu_\nu|^2}{4}\right)g^{n\bar n}
       +\frac{u_\nu^2}{4}\partial_\nu g^{n\bar n}.
\end{equation}

Along the fixed ray, the block determinant formula reads
$\det g=\det(g_{\alpha\bar\beta})_{\alpha,\beta<n}/g^{n\bar n}$.
At the boundary the tangential block is $(u_\nu/2)I$, so
differentiation of this identity and \eqref{eq:main} gives
$$
 \frac{\partial_\nu g^{n\bar n}}{g^{n\bar n}}
 =u_\nu+\frac2{u_\nu}
       \partial_\nu\sum_{\alpha<n}u_{\alpha\bar\alpha}.
$$
The polar-coordinate formula for the Euclidean Hessian gives
$$
 \sum_{\alpha<n}u_{\alpha\bar\alpha}(re_n)
 =\frac14\left\{
 r^{-2}\Delta_H[u(r\,\cdot)](e_n)+(2n-2)r^{-1}u_r(re_n)
 \right\},
$$
where $\Delta_H$ acts on the unit-sphere variable.
Using $u=0$ and $u_r=u_\nu$ at $r=1$, we obtain
\begin{align*}
 \partial_\nu\sum_{\alpha<n}u_{\alpha\bar\alpha}
 &=\frac14\Delta_Hu_\nu
   +\frac{n-1}{2}(u_{\nu\nu}-u_\nu)
 =\frac14\Delta_Hu_\nu
   +2(n-1)\left(u_{n\bar n}-\frac{u_\nu}{2}\right).
\end{align*}
Substitution in \eqref{eq:qnu-first}, using
$g^{n\bar n}=(u_\nu/2)^{n-1}/\mu$ at the boundary, yields
\begin{align*}
 \partial_\nu|\partial u|_g^2
 &=(n+1)u_\nu+\frac{(u_\nu/2)^{n-1}}{\mu}
 \left\{\frac{u_\nu^2}{2}\left(\frac{u_\nu}{2}-n\right)
 +\frac{u_\nu}{8}\Delta_Hu_\nu
 +\frac{n-1}{8}|\nabla_Hu_\nu|^2\right\}.
\end{align*}
The boundary identity in \eqref{eq:schur} also gives
$$
 \nabla_H|\partial u|_g^2
 =\frac{n+1}{2\mu}\left(\frac{u_\nu}{2}\right)^n\nabla_Hu_\nu.
$$
Combining the last two formulas, we find
\begin{align*}
 &\partial_\nu|\partial u|_g^2-
 \frac{\langle\nabla_Hu_\nu,\nabla_H|\partial u|_g^2\rangle}{2u_\nu}=(n+1)u_\nu+\frac{(u_\nu/2)^{n-1}}{\mu}
 \left\{\frac{u_\nu^2}{2}\left(\frac{u_\nu}{2}-n\right)
 +\frac{u_\nu}{8}\Delta_Hu_\nu
 -\frac14|\nabla_Hu_\nu|^2\right\}.
\end{align*}
On $\partial B$,
$$
 \partial_\nu e^{-u/(n+1)}=-\frac{u_\nu}{n+1},\qquad
 \nabla_He^{-u/(n+1)}=0,
$$
and the argument of the logarithm in \eqref{eq:log-positive}
equals $(u_\nu/2)^{n+1}/((n+1)\mu)$. Applying the product and
chain rules therefore gives
\begin{align*}
 &\partial_\nu\log\left[e^{-u/(n+1)}\left(1+\frac{|\partial u|_g^2}{n+1}\right)-1\right]-\frac1{2u_\nu}
       \left\langle\nabla_Hu_\nu,\nabla_H\log\left[e^{-u/(n+1)}\left(1+\frac{|\partial u|_g^2}{n+1}\right)-1\right]\right\rangle\\
 &\quad=\frac{2n\mu}{(u_\nu/2)^n}
       +2n\left(\frac{u_\nu}{2(n+1)}-1\right)
       +\frac{\Delta_Hu_\nu}{2u_\nu}
       -\frac{|\nabla_Hu_\nu|^2}{u_\nu^2}.
\end{align*}
By \eqref{eq:horizontal-ibp},
$$
 \int_{\partial B}\left(\frac{u_\nu}{2}\right)^{n-2}
       \Delta_Hu_\nu\,d\sigma_0
 =-\frac{n-2}{2}\int_{\partial B}
       \left(\frac{u_\nu}{2}\right)^{n-3}|\nabla_Hu_\nu|^2\,d\sigma_0.
$$
Substitution in \eqref{eq:flux} then yields
\begin{align*}
 &\int_B\Delta_g\log\left[e^{-u/(n+1)}\left(1+\frac{|\partial u|_g^2}{n+1}\right)-1\right]\det g\,dV_0\\
 &\quad=n^2\int_{\partial B}
 \left\{\frac{\mu}{u_\nu/2}-
 \left(\frac{u_\nu}{2}\right)^{n-1}
 \left(1-\frac{u_\nu}{2(n+1)}\right)\right\}d\sigma_0-\frac{n^2}{16}\int_{\partial B}
 \left(\frac{u_\nu}{2}\right)^{n-3}|\nabla_Hu_\nu|^2\,d\sigma_0,
\end{align*}
which proves \eqref{eq:log-flux}. These integrals are well defined
also for $n=2$, since $u_\nu>0$ on the compact sphere.
\end{proof}

\section{Symmetry and its consequences}\label{sec:normalized}

With the understanding of results above, we are ready to prove Theorem~\ref{thm:main}.
\begin{proof}[Proof of Theorem~\ref{thm:main}]
The parameter bound follows from \eqref{eq:mass-pohozaev}.

Suppose first that $0<\mu<n^n/(n+1)$, and let $a\in(0,n)$ and
$c\in(n,n+1)$ be the two positive roots of $s^n(1-s/(n+1))=\mu$.
By Lemma~\ref{lem:lower-slope}, $u_\nu/2\ge a$. For every $s\ge a$,
$$
 \left\{s^n\left(1-\frac{s}{n+1}\right)-\mu\right\}
 \left(\frac1s-\frac1c\right)\ge0,
$$
with equality only at $s=a$ or $s=c$: both factors are nonnegative
on $[a,c]$ and both are negative on $(c,\infty)$.
The Pohozaev identity gives
\begin{align*}
 &\int_{\partial B}\frac{
 \left(\frac{u_\nu}{2}\right)^n\left(1-\frac{u_\nu}{2(n+1)}\right)-\mu}
 {u_\nu/2}\,d\sigma_0=\int_{\partial B}
 \left\{\left(\frac{u_\nu}{2}\right)^n
 \left(1-\frac{u_\nu}{2(n+1)}\right)-\mu\right\}
 \left(\frac2{u_\nu}-\frac1c\right)d\sigma_0\ge0.
\end{align*}
Both boundary integrals in \eqref{eq:log-flux} are nonnegative,
while its left side is nonnegative by Proposition~\ref{prop:P}.
Thus all the integrals in \eqref{eq:log-flux} vanish. The continuous
nonnegative integrand in the last display is therefore zero, so
$u_\nu/2$ takes values only in $\{a,c\}$. Since the sphere is
connected, $u_\nu$ is constant. If $\mu=n^n/(n+1)$,
\eqref{eq:mass-pohozaev} gives $u_\nu/2\equiv n$, and
\eqref{eq:log-flux} again has zero left side. In both cases,
assertion~(3) of Proposition~\ref{prop:P} yields
$\nabla^{2,0}u=\frac1{n+1}\partial u\otimes\partial u$.

Fix $\zeta\in\partial B$. Since $u_\nu$ is constant there,
\eqref{eq:mass-pohozaev} gives
$$
 \mu=\left(\frac{u_\nu(\zeta)}{2}\right)^n
 \left(1-\frac{u_\nu(\zeta)}{2(n+1)}\right)>0,
 \qquad 0<\frac{u_\nu(\zeta)}2<n+1.
$$
The tensor identity implies
$$
 \nabla^{2,0}e^{-u/(n+1)}
 =-\frac{e^{-u/(n+1)}}{n+1}\nabla^{2,0}u
 +\frac{e^{-u/(n+1)}}{(n+1)^2}\partial u\otimes\partial u=0.
$$
Since the exponential is real, its $(0,2)$ covariant Hessian also
vanishes. For a real function $h$, the $(1,0)$ gradient
$\nabla_g^{1,0}h$ has components $\sum_jg^{i\bar j}h_{\bar j}$.
Metric compatibility gives
$$
 \partial_{\bar k}\left(\sum_jg^{i\bar j}
       \partial_{\bar j}e^{-u/(n+1)}\right)
 =\sum_jg^{i\bar j}\nabla_{\bar k}\nabla_{\bar j}e^{-u/(n+1)}=0.
$$
Thus $\nabla_g^{1,0}e^{-u/(n+1)}$ is holomorphic in $B$ and
smooth on $\overline B$.
On $\partial B$,
$\partial_{\bar j}e^{-u/(n+1)}=-u_\nu(\zeta)z_j/(2(n+1))$,
and the constant normal derivative makes the mixed block of $g$
vanish. Equation \eqref{eq:schur} therefore gives
$$
 \nabla_g^{1,0}e^{-u/(n+1)}
 =-\frac{(u_\nu(\zeta)/2)^n}{(n+1)\mu}
       \sum_i z_i\partial_{z_i}
 =-\frac1{n+1-u_\nu(\zeta)/2}\sum_i z_i\partial_{z_i}
 \quad\text{on }\partial B.
$$
The components of the difference between these holomorphic vector
fields vanish on $\partial B$. The maximum modulus principle
therefore extends the last equality throughout $B$.

The gradient identity gives
$$
 \bar\partial e^{-u/(n+1)}
 =-\frac{\bar\partial\left(\sum_i z_iu_i\right)}
 {n+1-u_\nu(\zeta)/2}.
$$
Since $\sum_i z_iu_i=u_\nu(\zeta)/2$ and $e^{-u/(n+1)}=1$ on
$\partial B$, the holomorphic function
$$
 \sum_i z_iu_i+\left(n+1-\frac{u_\nu(\zeta)}2\right)e^{-u/(n+1)}
$$
has constant boundary value $n+1$. Hence
$$
 \sum_i z_iu_i=n+1-
       \left(n+1-\frac{u_\nu(\zeta)}2\right)e^{-u/(n+1)}.
$$
The right side is real, so the real Euler derivative of $u$ equals
$2\sum_i z_iu_i$. Along every ray this yields
$$
 r\partial_r e^{u/(n+1)}=2e^{u/(n+1)}-
       \frac{2(n+1-u_\nu(\zeta)/2)}{n+1},
 \qquad e^{u/(n+1)}\big|_{r=1}=1.
$$
Solving on $0<r\le1$ gives, on every ray,
$$
 e^{u/(n+1)}=
 \frac{n+1-u_\nu(\zeta)/2+(u_\nu(\zeta)/2)r^2}{n+1}.
$$
Continuity includes the origin. Taking
$b=u_\nu(\zeta)/(2(n+1)-u_\nu(\zeta))>0$ proves \eqref{eq:radial}.

Conversely, for every $b>0$, the function in \eqref{eq:radial}
has tangential and radial complex Hessian eigenvalues
$(n+1)b/(1+b|z|^2)$ and $(n+1)b/(1+b|z|^2)^2$.
It therefore solves \eqref{eq:main} with the stated coefficient.
That coefficient increases on $0<b<n$, decreases on $b>n$,
and tends to zero at both ends. The asserted solution counts follow from this and Lemma \ref{lem:positive-roots}.
\end{proof}

For $u=\gamma\varphi$, the total mass $\gamma^n$ determines $b$
through $\gamma=(n+1)b/(1+b)$. For the alternative convention
$(i/\pi)\partial\bar\partial$, halving the potential and doubling
the coupling parameter preserve both the Monge--Amp\`ere measure
and the exponential density. The corresponding admissible range
of coupling parameters is $(0,2(n+1))$.

\begin{proof}[Proof of Corollary~\ref{cor:mean-field}]
Set $u=\gamma\varphi$. Then $u$ solves \eqref{eq:main} with
$\mu=\gamma^n/\int_Be^{-u}\,dV_0$ and has total
Monge--Amp\`ere mass $\int_B(dd^cu)^n=\gamma^n$.
By Theorem~\ref{thm:main}, $u_\nu/2$ is constant on $\partial B$
and belongs to $(0,n+1)$. The mass identity in
\eqref{eq:mass-pohozaev} gives $(u_\nu/2)^n=\gamma^n$ there,
hence $u_\nu/2=\gamma$. The formula in Corollary~\ref{cor:mean-field} follows from
\eqref{eq:radial}, and necessarily $0<\gamma<n+1$.

Conversely, for $0<\gamma<n+1$, the function in
Corollary~\ref{cor:mean-field} is
smooth and strictly plurisubharmonic on $\overline B$.
Direct differentiation and the mass identity in
\eqref{eq:mass-pohozaev}, applied to $u=\gamma\varphi$ with
$u_\nu/2=\gamma$, give
$$
 \det(\varphi_{i\bar j})
 =\left(1-\frac\gamma{n+1}\right)e^{-\gamma\varphi},
 \qquad
 \int_Be^{-\gamma\varphi}\,dV_0
 =\left(1-\frac\gamma{n+1}\right)^{-1}.
$$
These identities prove the stated equation.
\end{proof}

\begin{proof}[Proof of Corollary~\ref{cor:moser-trudinger}]
For $0<\gamma<n+1$, let $\varphi_\gamma$ be a maximizer of $G_\gamma$
over $\mathcal E^1(B)$ that is smooth up to the boundary and solves
\eqref{eq:mean-field}, as provided by \cite[Theorem~7.2 and Remark~7.3]{BB2022}.
Its complex Hessian is positive definite by plurisubharmonicity and
the positive density in \eqref{eq:mean-field}. Corollary~\ref{cor:mean-field}
identifies $\varphi_\gamma$ with the explicit radial solution.

Put $t=\gamma/(n+1)$ and $s=|z|^2$. The pushforward of $dV_0$
under $z\mapsto |z|^2$ is $ns^{n-1}\,ds$ on $[0,1]$.
By the preceding proof, $\int_B e^{-\gamma\varphi_\gamma}\,dV_0=(1-t)^{-1}$,
and the pushforward of $(dd^c\varphi_\gamma)^n$ under the same map is
$n(1-t)s^{n-1}(1-t+ts)^{-(n+1)}\,ds$. The substitution
$x=s/(1-t+ts)$ carries this measure to $n x^{n-1}\,dx$ on $[0,1]$
and gives $1-t+ts=(1-t)/(1-tx)$. Since
$t\varphi_\gamma=\log(1-t+ts)$, we obtain
\begin{equation}\label{eq:radial-functional}
 \gamma G_\gamma(\varphi_\gamma)
 =-n\int_0^1 x^{n-1}\log(1-tx)\,dx.
\end{equation}
Monotone convergence as $t\uparrow1$ and Tonelli's theorem give
\begin{align*}
 \lim_{\gamma\uparrow n+1}\gamma G_\gamma(\varphi_\gamma)
 &= n\int_0^1 x^{n-1}\int_0^x\frac{dy}{1-y}\,dx = \int_0^1\frac{1-y^n}{1-y}\,dy
  = \sum_{j=1}^n\frac1j.
\end{align*}
For a fixed $\varphi\in\mathcal E^1(B)$, maximality gives
$\gamma G_\gamma(\varphi)\le\gamma G_\gamma(\varphi_\gamma)$.
The defining approximants of $\varphi$ are nonpositive by the maximum
principle, so $\varphi\le0$. As $\gamma\uparrow n+1$,
$e^{-\gamma\varphi}$ increases to $e^{-(n+1)\varphi}$ and the energy
of $\varphi$ is finite. Monotone convergence therefore gives
$(n+1)G_{n+1}(\varphi)\le\sum_{j=1}^n1/j$.
For $u\in\mathcal E^1(B)$, dividing its defining sequence
by $n+1$ shows that $\varphi=u/(n+1)$ also belongs to $\mathcal E^1(B)$.
Homogeneity on the approximants and passage to their weak limits give
$(dd^c\varphi)^n=(n+1)^{-n}(dd^cu)^n$, which yields the asserted inequality.

The plurisubharmonic functions $u=(n+1)\varphi_\gamma$ belong to
$\mathcal E_0(B)$ and hence to $\mathcal E^1(B)$ because they are
smooth on $\overline B$ and vanish on $\partial B$.
For these functions the difference between
the logarithmic integral and the energy term is
$(n+1)G_{n+1}(\varphi_\gamma)$. Since $dV_0$ is a probability measure,
H\"older's inequality gives $G_{n+1}(\varphi_\gamma)\ge G_\gamma(\varphi_\gamma)$.
Together with the upper bound just proved, this shows that the
difference tends to $\sum_{j=1}^n 1/j$.
By \eqref{eq:radial-functional},
$$
 t\int_B(-\varphi_\gamma)(dd^c\varphi_\gamma)^n
 =-\log(1-t)-\gamma G_\gamma(\varphi_\gamma)\longrightarrow+\infty.
$$
The energy integral of $u=(n+1)\varphi_\gamma$ therefore diverges.
With any smaller energy coefficient the difference is unbounded above.
\end{proof}

Before closing this section, we mention the following result about the exterior problem, which is proved in \cite{CM26}.
Suppose $n\ge2$, $\mu>0$, and $v\in C^\infty(\overline{\Omega_{\C}};\mathbb R)$ has positive definite complex Hessian in $\Omega_{\C}:=\C^n\setminus\overline{B_1}$. Assume that
$$
 \det(v_{i\bar j})=\mu e^{-v}\quad\text{in }\Omega_{\C},\qquad
 v=0\quad\text{on }\partial B_1,
$$
and that $v(z)-k\log|z|^2$ is bounded uniformly as $|z|\to\infty$. Then for some $b>0$,
$$
 v(z)=k\log\frac{1+b|z|^2}{1+b},
 \qquad \mu=\frac{k^n b^n}{(1+b)^k}.
$$
The number of solutions is two for $0<\mu<n^n/(n+1)$, one for $\mu=n^n/(n+1)$ and zero for $\mu>n^n/(n+1)$.

\vspace{0.5cm}

\textbf{Data Availability Statement}
Data sharing is not applicable to this article because no datasets were generated or analyzed during the current study.


\textbf{Conflict of interest} The authors declare that they have no conflict of interest.

\textbf{Code availability}
Not applicable.

\textbf{AI disclosure}
The authors used AI tools to assist with this work.
The authors take full responsibility for the mathematical content and conclusions.

\end{document}